\documentstyle[amstex,12 pt]{article} \hfuzz=1.4pt %
\newtheorem{Th}{Theorem}[section]
\newtheorem{Prop}[Th]{Proposition}%[section]
\newcommand{\CQFD}
{%
\mbox{}%
\nolinebreak%
\hfill%
\rule{2mm}{2mm}%
\medbreak%
\par%
}

\newcommand{\R}{{\bf R}}

\newcommand{\C}{{\bf C}}

\newcommand{\r}{\rho}

\def\Ho{\vbox{\offinterlineskip\hbox{\kern 3pt$\scriptstyle\circ$}
\kern
1pt\hbox{$H$}}}
\begin{document}
\title{Finite Type Monge-Amp\`ere Foliations}

\author{Morris KALKA and Giorgio PATRIZIO}
\date{}%\today
\footnotetext{Much of this work was done while Kalka was
visiting the University in Florence. He would like to express his
thanks for the hospitality shown to him there. G. Patrizio thanks  Tulane
University where he was guest while  part of this research was completed and
acknowledges the support of MIUR and of the GNSAGA (INdAM).}

\maketitle
%\thispagestyle{empty}
%\newpage\addtocounter{page}{1}
\begin{abstract}
For plurisubharmonic solutions of the complex homogeneous Monge-Amp\`ere equation  whose level sets are hypersurfaces of finite type, in dimension $2$, it is shown that the Monge-Amp\`ere foliation is defined even at points of higher degeneracy. The result is applied to provide a positive answer to a question of Burns  on homogeneous polynomials whose logarithm satisfies the complex Monge-Amp\`ere equation and  to generalize the work of P.M. Wong  on the classification of complete weighted circular domains.\\

\end{abstract}
\section {Introduction}

Let $M$ be a Stein manifold of dimension $n$ and $\rho \colon M \to \R$ a
continuous plurisubharmonic exhaustion of class $C^\infty$ on $\{\rho >0\}$
with $d\rho \not= 0$ on $\{\rho >0\}$.  Furthermore suppose that
on $\{\rho >0\}$  the function $u=\log \rho$ satisfies the complex
homogeneous Monge-Amp\`ere equation,
\begin{align}(dd^c u)^n = 0.
\end{align}

On the open set, $P$, where $dd^c\rho >0$ the complex gradient  $Z$ of
$\rho$ (with respect to the K\"ahler metric with potential \(\rho\)) is
defined. In local coordinates $Z$ is given by
\begin{align}
Z=\sum_{\mu,\nu}\rho^{\mu\bar{\nu}}\rho_{\bar\nu}\frac{\partial}{\partial
z^{\mu}}
 \end{align}
where lower indeces denote derivatives and
$(\rho^{\mu\bar{\nu}})$ is defined by 
the relation $$\sum_{\nu=1}^{n} \rho^{\mu\bar{\nu}}\rho_{\alpha\bar{\nu}}=
\delta_{\alpha \mu}.$$

In \cite{Kalka-Patrizio}, it is shown that in case the function \(\r\) is
a real,
homogeneous polynomial on \(\C^n\) with convex level sets, the vector
field \(Z\) extends to all
of \(\C^n\) and yields a holomorphic foliation by parabolic Riemann surfaces.
Under the restrictive assumption of convexity we were able to give an answer to
a question of Burns (\cite{Burns}) and show that a plurisubhamonic, positive
homogeneous polynomial, \(\rho\) on  \(\C^n\) where
\(u=\log \rho\) satisfies \((dd^c u)^n = 0\) must be a \((p,p)\) polynomial. The
assumption of convexity allowed us to use the existence of complex geodesics for
the Kobayashi metric to extend the Monge-Amp\`ere foliation accross the weakly
pseudoconvex points and obtain a holomorphic foliation of \(\C^n\) by parabolic
Riemann surfaces. In the general pseudoconvex case the results about complex
geodesics for the Kobayashi metric are  not available. 
%\eject
 Here we make the
assumption that the level sets of our exhaustion are of finite type in order to
extend  the bundle \({\rm Ann}(dd^cu)\)  as an integrable complex line bundle on
\(M\). We are able to prove the existence of the extension in   dimension
\(n=2\). In this case we  show
that the extended bundle  yields a holomorphic foliation on \(\{\rho > 0\}\)
and, under the assumption that \(\{\rho > 0\}\) is connected, that the vector
field \(Z\) extends as a holomorphic vector field and that the zero set of
\(\rho\) is a point. Possibly under  more restrictive assumptions, we
believe that, the result holds for general \(n\). Our methods, however, use formulas 
that are only valid in dimension $2$, specifically the very simple 
form of the inverse of a $2\times 2$ matrix. We believe that there should be a more natural approach 
that avoids these $2-$dimensional calculations. 

Our extension and holomorphicity results allow us to answer Burns's question in
dimension \(2\) for any plurisuharmonic homogeneous polynomial. We are also
able to extend the results of P.M. Wong (\cite{Wong}) about classifying
complete weighted circular domains to the finite type case, for dimension
\(n=2\) assuming that the set \(\{\rho > 0\}\) is connected.

\section {Existence of the foliation}

We start recalling the well known notion of hypersurface of finite type in the
sense of Kohn \cite{Kohn} and Bloom-Graham \cite{BloomGraham}.
Let $S$ be a (real) hypersurface in a complex manifold $M$ and $p\in S$. Let $\rho
\colon U \to \R$ be a local defining function of $S$ at $p$, i.e.
$p\in U\cap S =\{z\in U \mid \rho(z)=0\}$ and $d\rho\neq 0$ on $U\cap S$.
Let  ${\cal L}^k$ denote the module over $C^{\infty}(U)$ functions generated by
all vector fields $W$ on $U$ in the holomorphic tangent bundle of $S$ (i.e.
with $\partial \rho(W)=0$ on $S$), their conjugates and their  brackets of length
at most $k$.  Then ${\cal L}^0$ is the module of vector fields on $U$ spanned
by the vector fields on $U$ holomorphically tangent to $S$ and their conjugates
and, for $k>0$,  ${\cal L}^k$ is the
the module of vector fields on $U$ spanned by brackets $[V,W]$ with
$V\in {\cal L}^{k-1}$ and $W\in{\cal L}^0$. The point $p\in S$ is  of {\sl type
m}  if $\partial \rho(W)=0$ for all $k=0,\dots, m-1$ and $W\in {\cal L}^k$
and there exists $Y\in {\cal L}^m$ with $\partial \rho(Y)\neq 0$. We say
that $S$ is of {\sl finite type} if for every $p\in S$ it is of type $m$ at
$p$, for some $m$ which depends on $p$.
\\ \\
\\
Throughout this section we shall make the following assumptions:
\\ \\
$(i)$ $M$ is a Stein manifold of dimension $2$;
\\ \\
$(ii)$ $\rho \colon M \to \R$ is
continuous plurisubharmonic exhaustion of class $C^\infty$ on $\{\rho >0\}$
with $d\rho \not= 0$ on $\{\rho >0\}$;
\\ \\
$(iii)$ $\rho$ is {\sl of finite type} i.e. all the level hypersurfaces $\{\rho
=c\}$,
for $c>0$ are of finite type (in the sense of Kohn \cite{Kohn} and Bloom-Graham
\cite{BloomGraham});
\\ \\
$(iv)$ the open set  $M_{*}=\{\rho >0\}$  is connected
\\ \\
$(v)$ on $\{\rho >0\}$  the function $u=\log \rho$ satisfies the complex
homogeneous Monge-Amp\`ere equation,
\begin{align}(dd^c u)^2 = 0.
\end{align}
\\
\textbf{Remark.} If the exhaustion $\rho$ is real analytic on
$M_{*}=\{\rho >0\}$, then $(iii)$ holds. In fact, since the level sets of the exhaustion $\rho$ are compact,  all level hypersurfaces $\{\rho =c\}$, for $c>0$ are of finite type  if $\rho$ is real analytic. Otherwise if $p$ were a point of  the real analytic hypersurface 
$S=\{\rho =c_{0}\}$ which is not of finite type, there  would   be a complex  variety  of positive dimension through $p$ which lies on S. This cannot happen is $S$ if compact (see   \cite{DiederichFornaess})
\\ \\
Suppose the assumptions $(i)$ to $(v)$ hold for a complex manifold $M$. On the
open set $P\subset M$ where $dd^c\rho >0$, the complex gradient  $Z$ of
$\rho$ (with respect to the K\"ahler metric with potential \(\rho\)) is
defined. In local coordinates $Z$ is given by
\begin{align}
Z=\rho^{\mu\bar{\nu}}\rho_{\bar\nu}\frac{\partial}{\partial
z^{\mu}}.
 \end{align}
It is well known (for instance see  \cite{Burns} or  \cite{Wong}) that, wherever it is defined, the vector field $Z$ is a non vanishing section of the line bundle \({\rm Ann}(dd^cu)\) and the Monge-Amp\`ere equation $(3)$ is
equivalent to 
\begin{align} \label{Zrho}
Z(\rho)=\rho.
\end{align}

We shall start showing that \({\rm Ann}(dd^cu)\) extends as a complex line
bundle on \(M\) and it is an integrable distribution spanned by a suitable
extension of $Z$.
\par

In any coordinate system  on an open set $U\subset \{\rho >0\}$ the local vector
field $L$ defined by
\begin{align}
\label{defL}
L= \rho_2 \partial_1 - \rho_1 \partial_2,
\end{align}
where subscripts for functions denote derivatives,
spans the holomorphic tangent space to level sets of  $\rho$.

The module ${\cal L}^k$  defined above is then  the module over $C^{\infty}(U)$
functions generated by $L$ and $\bar L$ and their  brackets of length at most
$k$.  Since the level sets of $\rho$ are hypersurfaces of finite
type, if $p\in U$, then there exists $k$ such that ${\cal L}^k_{|p}$ is the full complex
tangent space $T^{\C}_{p}(M)$.

On $P$, as $\partial\rho (Z)=Z(\rho)=\rho$, the complex gradient of $\rho$ is
transverse to holomorphic tangent bundle of any level set of $\rho$.
We shall extend  ${\rm Ann}(dd^cu)$  by suitably
choosing at the weakly pseudoconvex  points the ``missing direction''
recovered by the finite type assumption.
The following calculations, which will be useful later,  suggest also
that this is a reasonable approach.

\begin{Prop} If $U\subset \{\rho >0\}$ is an open coordinate set and
$L$ is defined by $(\ref{defL})$, on $P\cap U$ we have
\begin{align} \label{bracketL}
[L,\bar L]=&D(Z-\bar Z)\\
\label{bracket L,Z}
[L,Z]=&\phi_1L\\
[L,\bar Z]=&\psi_1L+\psi_2\bar L \label{bracket L,barZ}\\
[Z,\bar Z]=&\eta_1L+\eta_2\bar L\label{bracketZ}
\end{align}
where $D= \det
(\rho_{\mu\bar{\nu}})$ and  \(\phi_1,\phi_2,\psi_2\)
are  functions of class $C^{\infty}$.
\end{Prop}
{\bf Proof} The proof of (\ref{bracketL}) consists of a straightforward
calculation of
the bracket, together with the standard formula for the inverse of a
two by two matrix.
\begin{align*}
[L,\bar L]&= [\rho_2 \partial_1 - \rho_1 \partial_2,\rho_{\bar 2}
\partial_{\bar 1} - \rho_{\bar 1} \partial_{\bar
2}]\\
&=\left(\r_{\bar 1}\r_{2\bar 2}-\r_{\bar 2}\r_{2\bar 1}\right)\partial_1 +
\left(\r_{\bar 2}\r_{1\bar 1}-\r_{\bar 1}\r_{1\bar 2}\right)\partial_2 +
\\
&\hskip2cm\left( \r_2\r_{1\bar 2}-\r_1\r_{2\bar 2}\right)\partial_{\bar1} +
\left(\r_1\r_{2\bar 1}-\r_2\r_{1\bar 1}\right)\partial_{\bar2} \\
&=D(Z-\bar Z).
\end{align*}
In going to the last line we have used the identity: %correction%
\[
\begin{pmatrix} 
\r^{ 1 \bar1}& \r^{1 \bar2}\\
\r^{ 2 \bar1}& \r^{ 2 \bar2}
\end{pmatrix}
=\frac{1}{D}\begin{pmatrix}
\r_{2\bar 2} & -\r_{1\bar 2}\\
-\r_{2\bar 1} & \r_{1\bar 1}
\end{pmatrix}.
\]
\\ As for (\ref{bracket L,Z}), observe that
using the Monge-Amp\`ere equation \(Z(\r)=\bar Z(\r)=\r\)
and \(L(\r)=0\),
we have \([L,Z](\r)=LZ(\r)-ZL(\r)=0\). Thus,the \((1,0)\) vector
field \([L,Z]\) is tangent to each level set \(\{\r=c>0\}\) and
must therefore be a multiple of \(L\). The calculations of the
brackets \([L,\bar Z]\) and  \([Z,\bar Z]\) are entirely similar: both are
tangent to the level sets, so must be  linear combination of \(L\) and \(\bar
L\).

\CQFD

We can now prove

\begin{Th} The complex gradient
$Z$ defined on $P$ extends to a  non zero smooth vector field
on $\{\rho >0\}$ . The extension of the vector field $Z$  generates an
integrable complex line bundle ${\cal A}$ on $\{\rho >0\}$. Also the
Monge-Amp\`ere foliation defined by  ${\rm Ann}(dd^c u)$ extends to a foliation
defined by the distribution ${\cal A}$.
 \end{Th}
{\bf Proof} Let $p\in \{\rho >0\}\setminus P$ and suppose that
$L$ is the vector field defined by (\ref{defL}) on an open coordinate set
$U\subset\{\rho >0\}$ containing $p$.  Since  the boundary of each sublevel set
of $\rho$ is of finite type,  for some positive integer $m$ there exists a
$C^{\infty}$ nonzero vector field $Y\in {\cal L}^{m}$ with $\partial \rho(Y)\neq
0$ at $p$ and hence in a neighborhood of $p$ (which we may assume is $U$).
Since ${\cal L}^m$ is generated by $L$ and $\bar L$ and their  brackets of
length at most $m$,  we can
use (\ref{bracketL}),(\ref{bracket L,Z}), (\ref{bracket L,barZ}),
(\ref{bracketZ}) to  compute $Y$ on $U\cap P$ and conclude
\begin{align}\label{defY}
Y=\phi(Z-\bar Z)+ AL+B\bar L
\end{align}
for suitable functions $\phi, A, B$ of class $C^{\infty}$ on $U\cap P$.
If we define
\begin{align}\label{defV}
V=\frac{1}{2}(Y-iJY),
\end{align}
then on $U\cap P$ we have
\begin{align}\label{computeV}
V=\frac{1}{2}(Y-iJY)=\phi Z + AL.
\end{align}
We need now to compute and study the functions $\phi$ and
$A$. On  $U$ the function $\partial \rho(V)$ is smooth and non zero.
On the other hand on $U\cap P$, using (\ref{Zrho}) and (\ref{computeV})
\begin{align}
\partial \rho(V)=\phi\partial \rho(\rho)=\phi\rho
\end{align}
so that
\begin{align}
\phi=\frac{\partial \rho(V)}{\rho}
\end{align}
extends as a non zero function of class $C^{\infty}$ on all $U$.
\\
If, as before,  $D= \det (\rho_{\mu\bar{\nu}})$ denotes the determinant of the Levi
form of $\rho$, we define on $U\cap P$ the form
\begin{align}\label{defineOmega}
\Omega=\frac{\rho dd^{c}u}{D}+\frac{d\rho\wedge d^{c}\rho}{\rho}.
\end{align}
The form $\Omega$ itself does not extend on all $U$ but it is a sort of a
``singular'' metric which allows to compute and study the function $A$.
First of all,
using the equality
\begin{align}
dd^{c}u=\frac{dd^{c}\rho}{\rho}-
\frac{d\rho\wedge d^{c}\rho}{\rho^{2}},
\end{align}
observe that on $U\cap P$ 
\begin{align*}
\Omega^{2}&=({\rho dd^{c}u\over D}+{d\rho\wedge d^{c}\rho\over \rho})^{2}= {2\over D}dd^{c}u\wedge d\rho\wedge d^{c}\rho\\ &=
{2\over D}({ dd^{c}\rho\over \rho}- 
{d\rho\wedge d^{c}\rho\over \rho^{2}})\wedge d \rho\wedge d^{c}\rho
={2\over D}{dd^{c}\rho\wedge d \rho\wedge d^{c}\rho
\over \rho}.
\end{align*}
On the other hand, as $u$ satisfies the Monge-Amp\`ere equation,
we have
$$0=\rho^{3}(dd^{c}u)^{2}=\rho (dd^{c}\rho)^{2}-
2dd^{c}\rho\wedge d \rho\wedge d^{c}\rho$$
so that
$$\rho (dd^{c}\rho)^{2}=
2dd^{c}\rho\wedge d \rho\wedge d^{c}\rho.$$
We may conclude: %correction%
\begin{align*}\label{computeOmega2}
\Omega^{2}&={2\over D}{dd^{c}\rho\wedge d \rho\wedge d^{c}\rho
\over \rho}=
{\rho (dd^{c}\rho)^{2} \over D \rho}=
{(dd^{c}\rho)^{2} \over D}\\
&= \left({i\over 2\pi}\right)^{2}
dz_{1}\wedge dz_{2}\wedge d\bar z_{1}\wedge
d\bar z_{2}
\end{align*}
and therefore $\Omega^{2}$ is defined and of class $C^{\infty}$ on all $U$.
On  $U\cap P$ we have also
\begin{align}\label{OmegaZbarZ}
\Omega(Z,\bar Z)=\frac{\rho dd^{c}u(Z,\bar Z)}{D}+
\frac{d\rho\wedge d^{c}\rho(Z,\bar Z)}{\rho}=
\frac{\rho^{2}}{\rho}=\rho,
\end{align}
\begin{align}\label{OmegaZbarL}
\Omega(Z,\bar L)=\frac{\rho dd^{c}u(Z,\bar L)}{D}+
\frac{d\rho\wedge d^{c}\rho(Z,\bar L)}{\rho}=0,
\end{align}
where we have used the fact that $dd^{c}u(Z,\bar W)=0$ for any
vector field $W$ of type $(1,0)$ and that $d\rho\wedge d^{c}\rho(Z,\bar
L)=0$ holds, since $L\in \bf{Ker} \partial \rho$.    Using definitions and $Z(\rho)=\rho$,
a direct computation  gives
$dd^{c}\rho(L,\bar L)=D\rho$. Therefore,  using the fact that $L\in \bf{Ker} \partial \rho$, we have
\begin{align}\label{OmegaLbarL}
\Omega(L,\bar L)=\frac{\rho dd^{c}u(L,\bar L)}{D}+
\frac{d\rho\wedge d^{c}\rho(L,\bar L)}{\rho}=
\frac{\rho dd^{c}u(L,\bar L)}{D}=\rho.
\end{align}
The computations (\ref{OmegaZbarZ}),
(\ref{OmegaZbarL}), (\ref{OmegaLbarL}) show that $\Omega(Z,\bar Z)$,
$\Omega(Z,\bar L)$ and $\Omega(L,\bar L)$ extend of class $C^{\infty}$
throughout $U$. Thus $\Omega^{2}(V,\bar L,L,\bar L)$ is defined and of class $C^{\infty}$ throughout $U$ and on $U\cap P$ we have 
 following expression:
\begin{align*}\label{Omega2VbarV}
\Omega^{2}(V,\bar L,L,\bar L)&= \Omega^{2}(\phi Z+AL,\bar L,L,\bar L) = A \Omega^{2}(L,\bar L,L,\bar L)\\ &=
A[\Omega(L,\overline L)]^{2}= A\rho^{2}
\end{align*}
from which it follows that $A$ extends 
as   a $C^{\infty}$ functions  throughout $U$. 
\par
Therefore the complex gradient
$Z$ extends as a  $C^{\infty}$  vector field at the weakly pseudoconvex point
$p$ by setting $Z=\frac{1}{\phi}(V-AL)$.
The rest of the statement is now obvious.
\CQFD

\noindent
{\bf Definition}  We call the foliation defined by  ${\cal A}$ the {\sl extended Monge-Amp\`ere foliation}. 
\\
\\
\noindent
{\bf Remark} Let  $\ell$ be a leaf of the extended Monge-Amp\`ere foliation. At each point  $q\in \ell$ the tangent space $T_{q}\ell$ is the complex subspace of $T_{q}M$ spanned by  $Z(q)$ and hence
the leaf $\ell$ is a Riemann surface. The restriction $u_{|\ell}$  of the  function $u=\log \rho$ to $\ell$ is harmonic.  To see this, note that for $q\in M_{*}$, if $\ell$ is the leaf through $q$  and $\zeta$ is a holomorphic coordinate along $\ell$ in a coordinate neighborhood 
of $q$, then for some smooth function $\psi$ one has 
$\frac{\partial}{\partial \zeta}=\psi Z$. Hence the claim is equivalent to 
\[
dd^{c}u(\frac{\partial}{\partial \zeta},
\frac{\partial}{\partial \bar \zeta})=
\vert \psi\vert^{2}dd^{c}u(Z,\bar Z)=0
\]
where equality $dd^{c}u(Z,\bar Z)=0$  holds on $M^{*}$ as it is obvious  at  the points of the dense set $P$ where $dd^{c}\rho>0$.
\par
 Furthermore $u_{|\ell}$ is unbounded  above on $\ell$. Suppose  in fact that $u_{|\ell}$ were bounded above and let $r =\sup_{\ell} \rho>0$. Let $p\in \overline \ell$ with $\rho(p)=r$ and let $\ell_{p}$ be the leaf of the Monge-Amp\`ere foliation passing through $p$. Then necessarily $\ell_{p}\subset \{z\in M \mid \rho(z)\geq r\}$ otherwise the leaf $\ell$ would extend past $ \{z\in M \mid \rho(z)= r\}$. But then $p$ would be a local maximum for the harmonic function $u_{|\ell_{p}}$ which is impossible since $u_{|\ell_{p}}$ is not  constant.

 \begin{Th} The extended Monge-Amp\`ere foliation  on $M_{*}=\{\rho >0\}$ is
holomorphic.
\end{Th}
{\bf Proof} We have shown that the complex gradient $Z$ extends as a $C^{\infty}$-smooth vector field on
$M_{*}= \{z\in M \mid \rho(z)>0\}$. We denote again $Z$ this extension. The extended vector field  $Z$  is holomorphic along any leaf of the
Monge-Amp\`ere foliation. To see this, let us cover $M_{*}$ by coordinate
neighborhoods $U$ with $C^{\infty}$-smooth coordinates
$z_{1},z_{2}$ such that the intersection of a leaf of the Monge-Amp\`ere
foliation with $U$ is given by $z_{2}=c_{2}$ for suitable constants $c_{2}$ so
that $z_{1}$ is a holomorphic coordinate along each leaf. 
Since from (\ref{Zrho}) it follows that the leaves of the Monge-Amp\`ere foliation are transverse to the level set of $\rho$ and hence of $u$, in this coordinates, we have $du\wedge dz_{2}\neq 0$ so that, if $du=u_{1}dz_{1}+u_{2}dz_{2}$ then necessarily, with respect these coordinates, on $U$ we have  $u_{1}\neq 0$. Furthemore since $\partial \over \partial
z_{1}$ is tangent to the Monge-Amp\`ere foliation,  on  $U'=U\cap \{z\in
M_{*}\mid dd^{c}\rho(z)>0\}$ we have $Z=\varphi {\partial \over \partial
z_{1}}$ for some $C^{\infty}$-smooth function
$\varphi$. On the other hand on $U'$ we have $\rho=Z(\rho)=
\varphi {\partial \rho \over \partial z_{1}}$ from which it follows that on $U'$
we get
\begin{align}\label{extendZ two}
Z={1\over u_{1}} {\partial \over \partial z_{1}}.
\end{align}
Thus (\ref{extendZ two}) provides an expression of the $C^{\infty}$-smooth extension of $Z$ on $U$.
Since $u$ is harmonic along the leaves, we have that also  the extended vector field $Z$  is holomorphic along the intersection of $U$ with any leaf. Since $M_{*}$ is covered by such open coordinate neighborhoods, we got the claim.
\par
By continuity notice that we have also  on all $M_{*}$:
\begin{align}\label{MArho}
Z(\rho)=\rho
\end{align}
\\
Let $X={1\over 2}(Z+ \bar Z)$  and $Y={1\over 2i}(Z- \bar Z)$ be respectively
the real and the immaginary part of the vector field $Z$. Then, from
(\ref{MArho}) we conclude that
\begin{align}\label{XY two}
X(\rho)={1\over 4}\rho\hskip1cm{\rm and}\hskip1cm Y(\rho)=0.
\end{align}
Since $\rho\colon M_{*}\to (0,+\infty)$ is proper and the level sets of $\rho$
are compact, it follows from (\ref{XY two}) that $X$ and $Y$ are complete i.e.
the flows $\phi$ and $\psi$ of $X$ and $Y$ respectively, are both defined on
$\R$.
Fixed any point  $p\in M_{*}$, if $l_{p}$ is the leaf through $p$ then the map
$f \colon \C \to l_{p}$ defined by
\[f(t+is)=\phi(t,(\psi(s,p))\]
is holomorphic since $f'(t+is)=Z(f(t+is))$ and $Z$ is holomorphic along the leaf
and $f$ is non degenerated as $Z\neq0$ on $M_{*}$. Therefore the leaf $l_{p}$ is
parabolic. Since this is true for all leaves, Burns's Theorem  $3.2$ in
\cite{Burns} applies and therefore the Monge-Amp\`ere foliation is holomorphic.

\CQFD

Under the assumptions assumptions $(i)$ to $(v)$ stated at the beginning of the section, the main conclusion is the following:

\begin{Th} The vector field $Z$  extends as a  holomorphic vector
field on $M$ and the minimal set $\{z\in M \mid \rho(z)=0\}$ reduces to a
point.
\end{Th}
{\bf Proof} Let $p\in M_{*}= \{z\in M \mid \rho(z)>0\}$.  Since the extended
Monge-Amp\`ere foliation is holomorphic, there exists a coordinate neighborhood
$U$ of $p$ with holomorphic coordinates
$z_{1},z_{2}$ such that the intersection of a leaf of the Monge-Amp\`ere
foliation with $U$ is given by $z_{2}=c_{2}$ for suitable constants $c_{2}$ and
$z_{1}$ is a holomorphic coordinate along each leaf. Then, as we have seen
above, with respect to these coordinates  we get
\begin{align}\label{extendZbis}
Z={1\over u_{1}} {\partial \over \partial z_{1}}.
\end{align}
Thus, to show that $Z$ is holomorphic on $U$, we need only to prove
that $u_1$ is holomorphic. This is a simple consequence of the the fact $u$ is
solution of the Monge-Amp\`ere
equation. In fact since the vector field $Z$ generates at every point the
distribution ${\cal A}={\rm Ann} (dd^{c}u)$, then for every vector field
$X=\sum X_{\alpha}{\partial\over\partial z_{\alpha}}$ of type $(1,0)$ on $U$, it
must be that $dd^{c}u(Z,\bar X)=0$. In local coordinates, using the expression
(\ref{extendZbis}), this is equivalent to
\[
u_{1\bar1}{1\over u_{1}}\overline{X^{1}}+
u_{1\bar2}{1\over u_{1}}\overline{X^{2}}=0
\]
which, since $X$ is arbitrary and $u_{1}\neq 0$, implies that
$u_{1\bar1}=u_{1\bar2}=0$  and therefore $u_{1}$, and hence $Z$, is holomorphic
on $U$. Since $p$ is arbitrary in $M_{*}$, then
$Z$ is holomorphic on $M_{*}$. Since $M$ is Stein, if $M_{*}$ is connected 
and $\{z\in M \mid \rho(z)=0\}$ is compact, then the vector field $Z$ extends
holomorphically to all $M$.
\par
Finally we like to show that the minimal set $\{z\in M \mid \rho(z)=0\}$ of
$\rho$ reduces to one point. Firstly observe that it cannot be empty as
otherwise $u=\log\rho$ would be a smooth exhaustion which solves the complex
homogeneous Monge-Amp\`ere equation on a Stein manifold. But such function does
not exist (Theorem 1.1 of \cite{lempert-szoke}). The holomorphic flow
$\Psi\colon \C\times M_{*}\to M_{*}$ on $M_{*}$ of the vector field $Z$ extends
to a holomorphic flow $\Psi\colon \C\times M\to M$ as $\{z\in M \mid
\rho(z)=0\}$ is compact. On the other hand for any $p\in \{z\in M \mid
\rho(z)=0\}$, since the flow of any point in $M_{*}$ is contained in $M_{*}$ we
must have
$\Psi(\C\times\{p\})\subset \{z\in M \mid \rho(z)=0\}$, which, by Liouville
theorem
implies that $\Psi(\C\times\{p\}=\{p\})$. Thus $Z(p)=0$. But then
$\{z\in M \mid \rho(z)=0\}=\{z\in M \mid Z(z)=0\}$ is a compact analytic set,
i.e. a finite set of points. To finish our proof it is enough to show that the
minimal set of $\rho$ is connected. This is achieved repeating verbatim the
argument
at page 357 of \cite{Burns}. We give here an outline of the proof.  Suppose that
for some compact, non empty disjoint compact subset $K_{1}, K_{2}\subset M$ we
have $\{z\in M \mid \rho(z)=0\}=K_{1}\cup K_{2}$. Let $V_{1}, V_{2}$ disjoint
open sets with $K_{i}\subset V_{i}$ for $i=1,2$.  For $r>0$ sufficiently small
we have $S_{i}=\rho^{-1}(r)\cap V_{i}\neq \emptyset$ for $i=1,2$.
If $G\colon \R\times M_{*}\to  M_{*}$ is the flow of the real part $X$ of the
vector field $Z$, then $G(\R\times\rho^{-1}(r))=M_{*}$. On the other hand
if $U_{i}=G(\R\times S_{i})$ for $i=1,2$, then $U_{1},U_{2}$ are open disjoint
subsets with $U_{1}\cup U_{2}=M_{*}$ contradicting the fact that $M_{*}$ is
connected.

\CQFD

\section {Burns's Problem}

Taking advantage of the results of section 2 and using the final arguments  of
\cite{Kalka-Patrizio}, we can provide a
complete solution to the problem proposed by  D. Burns \cite{Burns}
in dimension $2$. Namely we can show the following:

\begin{Th} Let $\r$ be a  positive homogeneous polynomial on
${\bf C}^2$  such that $u=\log \r$ is plurisubharmonic and satisfies
\begin{align} \label{MA}(dd^c u)^2 = 0.
\end{align}
Then  $\r$ is a homogeneous polynomial of bidegree $(k,k)$.
\end{Th}
{\bf Proof}
From the results of Section 2 it follows that the complex gradient $Z$, defined 
on $P$, using as usual  the summation convention, by
\begin{align}\label{zedbis}
Z=Z^{\mu}{\frac{\partial\,} {\partial z^{\mu}}} =\rho^{\mu\bar{\nu}}\rho_{\bar\nu}\frac{\partial}{\partial
z^{\mu}},
 \end{align}
extends holomorphically to all ${\bf C}^2$. On the other hand,  (\ref{zedbis})
shows that $Z$ is homogeneous of degree one on a dense subset of
${\bf C}^2$, so that $Z$ is   in fact linear.
\\
We can get an explicit expression  for $Z$ using a bidegree
argument. Suppose that
$$\r=\sum_{l+m=2k}\r^{l,m}$$
is the decomposition of $\r$ in homogeneous polynomial of bidegree $(l,m)$.
The Monge-Amp\`ere equation (\ref{MA}) is equivalent on $P$ to 
\[Z(\r)=\r^{\mu\bar \nu}\r_{\bar\nu}\r_{\mu}= \r=\sum_{l+m=2k}\r^{l,m}. \]
Differentiating  this equation with respect to $\bar z^\alpha$, we get
$$\sum_{\scriptstyle l+m=2k \atop \scriptstyle m\geq 1}\r^{l,m}_{\bar\alpha}
=\r_{\bar\alpha}=
\left(Z^{\mu}\r_{\mu}\right)_{\bar\alpha}=
\sum_{\scriptstyle l+m=2k\atop\scriptstyle  l,m\geq
1}Z^{\mu}\r^{l,m}_{\mu\bar\alpha}.$$
Comparing bidegrees, we conclude:
$$0=\r^{0,2k}_{\bar\alpha}=\r^{2k,0}_{\bar\alpha}$$
and, for every $l,m$ with $l+m=2k, l,m\geq 1$,
$$\r^{l,m}_{\bar\alpha}=Z^{\mu}\r^{l,m}_{\mu\bar\alpha}$$
where, as usual we use the summation convention.
Let $w=(w^1,\dots,w^n)$ be such that $dd^c\r_w>0$; then $dd^c\r^{k,k}_w>0$
(see for instance Lemma $2.1$ in \cite{Kalka-Patrizio}) and,
if $(\r^{k,k})^{\bar\alpha\mu}$ is defined by 
the relation 
$$\sum_{\nu=1}^{n} (\r^{k,k})^{\mu\bar{\nu}}(\r^{k,k})_{\alpha\bar{\nu}}=
\delta_{\alpha \mu},$$
then
\begin{align}\label{newzed} 
Z^{\mu}_w=
(\r^{k,k})^{\bar\alpha\mu}(w)\r^{k,k}_{\bar\alpha}(w)={1\over k}w^k
\end{align}
where the last equality is immediate from $(5)$ as $\r^{k,k}_{\bar\alpha}$
is homogeneous of degree $(k,k-1)$ and hence
$\r^{k,k}_{\mu\bar\alpha}(W)w^\mu=k\r^{k,k}_{\bar\alpha}$.
By continuty (\ref{newzed}) holds on all ${\bf C}^2$.
\par
As consequence the (extended) Monge-Amp\`ere foliation associated to
$u=\log\r$ is given by the foliation of  ${\bf C}^2\setminus\{0\}$ in lines
through the origin. It follows that the restriction of $\log \r$ to any complex
line through the origin is harmonic with a logarithmic singularity of weight
$2k$ at the origin
and therefore for any $0\not=z\in{\bf C}^2$  and $\lambda \in {\bf C}$
$$\log \r(\lambda z)=2k\log \vert \lambda \vert + O(1).$$
Hence the restriction of $\r$ to any complex line through the origin is
homogeneous of bidegree $(k,k)$ and therefore $\r$ is a homogeneous polynomial
of bidegree $(k,k)$ on
${\bf C}^2$.

\CQFD

\section {Classification of solutions}

Given what we have already proven, we see that the argument of P.M. Wong in
(\cite{Wong}) goes over and we  are able to prove  the following:
%correction%

\begin{Th} Let $M$ be a Stein manifold of dimension $2$ equipped with a
continuous plurisubharmonic exhaustion $\rho \colon M \to \R$  of class
$C^\infty$ on $\{\rho >0\}$  with $d\rho \not= 0$ on $\{\rho >0\}$ such that
\\
$(i)$ $\rho$ is {\sl of finite type} i.e. all the level hypersurfaces $\{\rho=
c\}$,
for $c>0$ are of finite type (in the sense of Kohn \cite{Kohn} and Bloom-Graham
\cite{BloomGraham}),
\\
$(ii)$  the open set  $M_{*}=\{\rho >0\}$  is connected,
\\
$(iii)$ on $\{\rho >0\}$  the function $u=\log \rho$ satisfies the complex
homogeneous Monge-Amp\`ere equation,
\begin{align} \label{MAbis} (dd^c u)^2 = 0.
\end{align}
Then there exists a biholomorphic map $\Phi\colon {\bf C}^2\to M$,
such that, for suitable positive real numbers $c_{1},c_{2}$, the pull back
$\r_{0}=\r\circ\Phi$ of the exahustion $\r$ satisfies
\begin{align}\label{rholinerized}
\rho_{0}(e^{c_{1}\lambda}z_{1},e^{c_{2}\lambda}z_{2})=
\vert e^{\lambda}\vert^{2}\rho_{0}(z).
\end{align}
so that the sublevelsets of $\rho$ are biholomorphic to generalized weighted
circular domains.
\end{Th}
{\bf Proof}
Because of the results of Section 2, our hypotheses imply that the
Monge-Amp\`ere foliation associated to $u=\log \rho$ extends to a holomorphic
foliation throughout $M_{*}=\{\rho >0\}$.
\par
Furthemore the complex gradient vector field $Z$, locally defined by
\begin{align}\label{zed}
Z=\rho^{\mu\bar{\nu}}\rho_{\bar\nu}\frac{\partial}{\partial
z^{\mu}},
 \end{align}
on the set where $dd^{c}\rho>0$, extends holomorphically to all $M$ so that the
extension, which we shall keep denoting $Z$, is tangent to the leaves of the
extended Monge-Amp\`ere foliation  on $M_{*}$.
\par
The equality $Z(\rho)=\rho$,
which is equivalent to the Monge-Amp\`ere equation (\ref{MAbis}) on the set
where $dd^{c}\rho>0$, holds, by continuity on all $M$.
\par
Finally we have that the minimal set of the function $\rho$ reduces to a point:
$\{\rho =0\}=\{O\}$. We shall call the point $O$ the {\sl center} of $M$.
The proof rests on the ground that the argument in
(\cite{Wong}) uses only these facts and no other consequence of strict
pseudoconvexity. We outline the main steps of the rest of the proof there to
underline the minor varitions needed under our weaker assumptions.
\\ \\
{\bf Step 1:} For any $r_{1}, r_{2}>0$ the level sets $M<r_{i}>=\{\rho=r_{i}\}$
for $i=1,2$ are $CR$ isomorphic and the sublevel sets $M(r_{i})=\{\rho<r_{i}\}$
for $i=1,2$ are biholomorphic.
\\ \\
Step 1 is consequence of a Morse Theory type of argument. The flow of the real
part $X=Z+\bar Z$ of the (extended) complex gradient vector field maps level
defines a local group of biholomorphims when $Z$ is holomorphic and maps level
sets of $\r$ into other level sets. Details of the argument can be found on
page 24 of (\cite{Patrizio-Wong}). There it was assumed   $\r$ to be strictly
plurisubharmonic  merely to ensure that the vector field $Z$ is defined
everywhere which in our case we prove by other means. The claim then follows
from a standard application of Bochner-Hartogs extension theorem.
\\ \\
{\bf Step 2:} There exists a coordinate neighborhood $U$ of the center $O$, with
coordinates $z_{1},z_{2}$ centered at  $O$, such that for suitable positive real
numbers $c_{1}, c_{2}>0$, on $U$, the vector field $Z$ has the following
expression:
\begin{align}\label{zedlinerized}
Z=c_{1}z_{1}{\partial \over \partial z_{1}}+
c_{2}z_{2}{\partial \over \partial z_{2}}. \end{align}
\\
Here we can be repeat the proof on page 248 of \cite{Wong} which goes as
follows: the equation $Z(\rho)=\rho $ implies that the one parameter group
associated with imaginary part of $Z$ preserves the level sets of $\rho$ and
therefore have compact closure in the automorphism group of each sublevel set.
A theorem in Bochner-Martin (\cite{Bochner-Martin}, Chapter III) then implies
the given form of vector field $Z$.
\\ \\
{\bf Step 3:} There exist  suitable positive real numbers $c_{1},c_{2}$ so that
for $z=(z_{1},z_{2})\in U$ and $\lambda \in \C$ so that
$(e^{c_{1}\lambda}z_{1},e^{c_{2}\lambda}z_{2})\in U$, we have:
\begin{align}\label{rholinerized}
\rho(e^{c_{1}\lambda}z_{1},e^{c_{2}\lambda}z_{2})=
\vert e^{\lambda}\vert^{2}\rho(z).
\end{align}
\\
The proof is as in page 249 of \cite{Wong} and consists of integrating the
vector field $Z$ explicitly. Again, one uses only the the fact that
$Z(\rho)=\rho$ in this integration.
\\ \\
What is left to prove is that there is a global biholomorphism of
${\bf C}^2$ onto $M$. This follows from the fact that we can define the
biholomorphism on a neighborhood of the origin to a neighborhood of the center
of $M$ and then flow in each manifold with the real part of $Z$. More
specifically, from Step 3 we may conclude for  $\epsilon>0$ small enough the
sublevel set $M(\epsilon)$ of $\rho$ is biholomorphic to a (fixed) weighted
circular domain
$$G(\epsilon)=\{z\in  {\bf C}^2 \mid \r_{0}(z)<\epsilon\}$$
for some plurisubharmonic exhaustion $\r_{0}$ of ${\bf C}^2$
which satisfies (\ref{rholinerized}).
If $\varphi : M(\epsilon)\to G(\epsilon)$ is the biholomorphism of Step 3, then
$\r_{0}=\r\circ \varphi$ on $G(\epsilon)$, it is defined on all ${\bf C}^2$ and
it satisfies $(i)$, $(ii)$, $(iii)$.
By Step 1 it follows that the flows of the real parts of the complex gradients
of $\r$ and $\r_{0}$ are biholomorhisms of sublevel sets
The required biholomorphic map $\Phi \colon {\bf C}^2 \to M$ can
be defined by composition of $\varphi$ and the flows of the real parts
of the complex gradients.
\CQFD

\end{document}